\documentclass[11pt]{amsart}
\usepackage{amsfonts}
\usepackage{amscd}
\usepackage{amsmath}
\usepackage{epsfig}

\newtheorem{definition}{Definition}
\newtheorem{theorem}{Theorem}

\newtheorem{corollary}{Corollary}
\newtheorem{remark}{Remark}
\newtheorem{lemma}{Lemma}

\def\be{\begin{equation}}
\def\ee{\end{equation}}
\def\ben{\begin{displaymath}}
\def\een{\end{displaymath}}
\def\baa{\begin{eqnarray}}
\def\eaa{\end{eqnarray}}

\def\ba{\begin{array}}
\def\ea{\end{array}}
\def\CP1{{\bf CP}^1}

\begin{document}

\title{Determinant of pseudo-laplacians}
\author[Tayeb Aissiou]{Tayeb Aissiou}
\email{aissiou@math.mcgill.ca}
\address{Department of Mathematics and Statistics, Concordia University\\
1455 de Maisonneuve Blvd. West \\
Montreal, Quebec H3G 1M8 Canada}

\author[L. Hillairet]
{Luc Hillairet}
\email{Luc.Hillairet@math.univ-nantes.fr}
\address{UMR CNRS 6629-Universit\'{e} de Nantes, 2 rue de la Houssini\`{e}re \\
BP 92 208, F-44 322 Nantes Cedex 3, France}

\author[A. Kokotov]
{Alexey Kokotov}
\email{alexey@mathstat.concordia.ca}
\address{Department of Mathematics and Statistics, Concordia University\\
1455 de Maisonneuve Blvd. West \\
Montreal, Quebec H3G 1M8 Canada}

\maketitle

 \vskip0.5cm {\bf Abstract.}
 We derive comparison formulas relating the zeta-regularized determinant of an arbitrary self-adjoint extension of the Laplace operator
 with domain $C^\infty_c(X\setminus \{P\})\subset L_2(X)$ to the zeta-regularized determinant of the Laplace operator on $X$.
 Here  $X$ is a compact Riemannian manifold of dimension $2$ or $3$; $P\in X$.
  \vskip0.5cm

\section{Introduction}

Let $X_d$ be a complete Riemannian manifold of dimension $d\geq 2$ and let $\Delta$ be the (positive) Laplace operator on $X_d$. Choose a point $P\in X_d$ and consider $\Delta$ as an unbounded symmetric operator in the space $L_2(X_d)$ with domain $C^\infty_c(X_d\setminus\{P\})$. It is well-known that thus obtained operator is essentially self-adjoint if and only if $d\geq 4$. In case $d=2, 3$ it has deficiency indices $(1, 1)$ and there exists a one-parameter family $\Delta_{\alpha, P}$ of its self-adjoint extensions (called pseudo-laplacians; see \cite{YCdV}). One of these extensions (the Friedrichs extension $\Delta_{0, P}$) coincides with the self-adjoint operator $\Delta$ on $X_d$. In case $X_d=R^d$, $d=2,3$ the scattering theory for the pair $(\Delta_{\alpha, P}, \Delta)$ was extensively studied in the literature (see e. g.,  \cite{Albeverio}). The spectral theory of the operator $\Delta_{\alpha, P}$ on  compact manifolds $X_d$ $(d=2, 3)$ was studied in \cite{YCdV}, notice also a recent paper \cite{Ueber} devoted to the case, where $X_d$ is a compact Riemann surface equipped with Poincar\'{e} metric.

The zeta-regularized determinant of Laplacian on a compact Riemannian manifold was introduced in \cite{RS} and since then was studied and used in an immense number of papers in string theory and geometric analysis, for our future purposes we mention here the memoir \cite{Fay}, where the determinant of Laplacian is studied as a functional on the space of smooth Riemannian metrics on a compact two-dimensional manifold, and the papers \cite{FG2003} and \cite{Kum}, where the reader may find explicit calculation of the determinant of  Laplacian  for  three-dimensional flat tori and  for the sphere $S^3$ (respectively).

The main result of the present paper is a comparison formula relating ${\rm det}(\Delta_{\alpha, P}-\lambda)$ to ${\rm det}(\Delta-\lambda)$,
for $\lambda \in {\mathbb C}\setminus \left({\rm Spectrum} (\Delta)\cup {\rm Spectrum}(\Delta_{\alpha, P})\right)$.

It should be mentioned that in case of two-dimensional manifold the zeta-regularization of  ${\rm det}(\Delta_{\alpha, P}-\lambda)$ is not that standard, since the corresponding operator zeta-function has logarithmic singularity at $0$.

It should be also mentioned that in the case when the manifold $X_d$ is flat in a vicinity of the point $P$ we deal with a very special case of the situation (Laplacian on a manifold with conical singularity) considered in \cite{Loya}, \cite{KLP1}, \cite{KLP} and, via other method, in  \cite{HK}. The general scheme of the present work is close to that of \cite{HK}, although some calculations from \cite{KLP} also appear very useful for us.

{\bf Acknowledgements.}The work of T. A. was supported by FQRNT. Research of A. K. was supported by NSERC.

\section{Pseudo-laplacians, Krein formula and scattering coefficient}
Let $X_d$ be a compact manifold of dimension $d=2$ or $d=3$; $P\in X_d$ and $\alpha\in [0, \pi)$. Following Colin de Verdi\`{e}re \cite{YCdV}, introduce the set
$${\mathcal D}(\Delta_{\alpha, P})=\{f\in H^2(X_d\setminus \{P\}): \exists c\in {\mathbb C}: {\text \ in\  a\  vicinity\  of\ } P {\text \ one \ has}$$
\begin{equation}\label{as} f(x)=c(\sin \alpha \cdot G_d(r)+\cos\alpha)+o(1) {\text \ as\ } r\to 0\}\,,\end{equation}
where
$$H^2(X_d\setminus \{P\})=\{f\in L_2(X_d): \exists C\in {\mathbb C}: \Delta f- C\delta_{P}\in L_2(X_d)\}\,,$$
$r$ is the geodesic distance between $x$ and $P$ and
$$G_d(r)=\begin{cases}
\frac{1}{2\pi}\log r, \ \  d=2\\
-\frac{1}{4\pi r}, \ \ d=3.
\end{cases}$$

Then (see \cite{YCdV}) the self-adjoint extensions of symmetric operator $\Delta$ with domain $C^\infty_c(X_d\setminus\{P\})$ are the operators $\Delta_{\alpha, P}$
with domains ${\mathcal D}(\Delta_{\alpha, P})$ acting via $u\mapsto \Delta u$. The extension $\Delta_{0, P}$ coincides with the Friedrichs extension and is nothing but the self-adjoint Laplacian on $X_d$.

Let $R(x, y; \lambda)$ be the resolvent kernel of the self-adjoint Laplacian $\Delta$ on $X_d$.

 Following \cite{YCdV} define the scattering coefficient $F(\lambda; P)$ via
\begin{equation}\label{coeffF} -R(x, P; \lambda)=G_d(r)+F(\lambda; P)+o(1)\end{equation}
as $x\to P$.
(Notice that  in \cite{YCdV} the resolvent is defined as $(\lambda-\Delta)^{-1}$, whereas for us it is $(\Delta-\lambda)^{-1}$. This results in the minus sign in (\ref{coeffF}).)

As it was already mentioned the deficiency indices of the symmetric operator $\Delta$ with domain $C^\infty_c(X_d\setminus\{P\})$ are $(1, 1)$, therefore, one has the following Krein formula (see, e. g., \cite{Albeverio}, p. 357) for the resolvent kernel, $R_\alpha (x, y; \lambda)$,  of the self-adjoint extension $\Delta_{\alpha, P}$:
\begin{equation}\label{Krein}
R_\alpha(x, y; \lambda)=R(x, y; \lambda)+k(\lambda; P)R(x, P; \lambda)R(P, y; \lambda)
\end{equation}
with some $k(\lambda; P)\in {\mathbb C}$.

The following Lemma relates $k(\lambda; P)$ to the scattering coefficient $F(\lambda; P)$.
\begin{lemma}\label{L1} One has the relation
\begin{equation}\label{rel}
k(\lambda; P)=\frac{\sin \alpha}{F(\lambda; P)\sin \alpha-\cos\alpha}\,.
\end{equation}
\end{lemma}
{\bf Proof.} Send $x\to P$ in (\ref{Krein}), observing that $R_\alpha (\,\cdot\, , y; \lambda)$ belongs to ${\mathcal D}_{\alpha, P}$,  make use of (\ref{as}) and (\ref{coeffF}), and then compare the coefficients near $G_d(r)$ and the constant terms in the asymptotical expansions at the left and at the right. $\square$

It follows in particular from the Krein formula that the difference of the resolvents $(\Delta_{\alpha, P}-\lambda)^{-1}-(\Delta-\lambda)^{-1}$
is a rank one operator. The following simple Lemma is the key observation of the present work.
\begin{lemma}
One has the relation
\begin{equation}\label{trace}
{\rm Tr}\,\left( (\Delta_{\alpha, P}-\lambda)^{-1}-(\Delta-\lambda)^{-1}  \right)=\frac {F_\lambda '(\lambda; P)\sin\alpha}{\cos\alpha-F(\lambda; P)\sin\alpha}\,.
\end{equation}
\end{lemma}

{\bf Proof.} One has
$$-F_\lambda'(\lambda; P)=\frac{\partial R(y, P; \lambda)}{\partial \lambda}\Big|_{y=P}=\lim_{\mu\to\lambda}\frac{R(y, P; \mu)-R(y, P; \lambda)}{\mu-\lambda}$$
Using resolvent identity we rewrite the last expression as
$$\lim_{\mu\to \lambda}\int_{X_d}R(y, z; \mu)R(P, z; \lambda)dz\Big|_{y=P}=\int_{X_d}[R(P, z; \lambda)]^2dz$$
From (\ref{Krein}) it follows that
$$[R(P, z; \lambda)]^2=\frac{1}{k(\lambda; P)}\left(R_{\alpha, P}(x, z; \lambda)-R(x, z; \lambda)\right)\Big|_{x=z}\,.$$
This implies
$$-F_\lambda'(\lambda; P)=\frac{1}{k(\lambda, P)}{\rm Tr}\,\left( (\Delta_{\alpha, P}-\lambda)^{-1}-(\Delta-\lambda)^{-1}  \right)$$
which, together with Lemma \ref{L1}, imply (\ref{trace}).$\square$

Introduce the domain
$$\Omega_{\alpha, P}={\mathbb C}\setminus \{\lambda-it, \lambda\in {\rm Spectrum}\,(\Delta)\cup{\rm Spectrum}\,(\Delta_{\alpha, P}); t\in (-\infty, 0]\}\,.$$
Then in $\Omega_{\alpha, P}$ one can introduce the function
\begin{equation}\label{xsi}\tilde\xi(\lambda)=-\frac{1}{2\pi i} \log(\cos\alpha-F(\lambda; P)\sin\alpha)\end{equation}
(It should be noted that the function $\xi=\Re (\tilde \xi)$ is the spectral shift function of $\Delta$ and $\Delta_{\alpha, P}$.)
One can rewrite (\ref{trace})
as
\begin{equation}\label{trace1}
{\rm Tr}\,\left( (\Delta_{\alpha, P}-\lambda)^{-1}-(\Delta-\lambda)^{-1}  \right)=2\pi i\tilde \xi'(\lambda)
\end{equation}
\section{Operator zeta-function of $\Delta_{\alpha, P}$}

Denote by $\zeta(s, A)$ the zeta-function
$$\zeta(s, A)=\sum_{\mu_k\in {\rm Spectrum}\,(A)}\frac{1}{\mu_k^s}$$
of the operator $A$. (We assume that the spectrum of $A$ is discrete and does not contain $0$.)

 Take any $\tilde \lambda$ from ${\mathbb C}\setminus ({\rm Spectrum}\,(\Delta_{\alpha, P})\cup {\rm Spectrum}\,(\Delta)))$. From the results of \cite{YCdV} it follows that the function $\zeta(s, \Delta_{\alpha, P}-\tilde \lambda)$ is defined for sufficiently large $\Re s$. It is well-known that $\zeta(s, \Delta-\tilde \lambda)$ is meromorphic in ${\mathbb C}$.

The proof of the following lemma coincides verbatim with the proof of Proposition 5.9 from \cite{HK}.
\begin{lemma}\label{lem}
Suppose that the function $\tilde \xi'(\lambda)$ from (\ref{trace1}) is $O(|\lambda|^{-1})$ as $\lambda\to -\infty$.
Let $-C$ be a sufficiently large negative number and let $c_{\tilde \lambda,\epsilon}$ be a contour encircling the cut
  $c_{\tilde \lambda}$  which starts from $-\infty+0i$, follows the real line till $-C$ and then goes to $\tilde \lambda$ remaining in $\Omega_{\alpha, P}$. Assume that
  ${\rm dist}\,(z, c_{\tilde \lambda})=\epsilon$ for any $z\in c_{\tilde \lambda,\epsilon}$.
Let also
$$\zeta_2(s)=\int_{c_{\tilde \lambda, \epsilon; 2}}(\lambda-\tilde \lambda)^{-s}\tilde\xi'(\lambda)d\lambda,$$
where the the integral at the right hand side is taken over the part $c_{\tilde \lambda, \epsilon; 2}$ of the contour $c_{\tilde \lambda,\epsilon}$ lying in the half-plane $\{\lambda: \Re\lambda>-C\}$.
Let $$\hat{\zeta}_2(s)=\lim_{\epsilon\to 0}\zeta_2(s)=2i\sin(\pi s)\int_{-C}^{\tilde\lambda}(\lambda-\tilde\lambda)_0^{-s}\tilde \xi'(\lambda)\,d\lambda\,,$$
where $(\lambda-\tilde\lambda)_0^{-s}=e^{-i\pi s}\lim_{\lambda\downarrow c_{\tilde \lambda}}(\lambda-\tilde\lambda)^{-s}$.
Then the function

\begin{equation}\label{MAIN}
R(s, \tilde \lambda)= \zeta(s, \Delta_{\alpha, P}-\tilde\lambda))-\zeta(s, \Delta-\tilde\lambda)-2i\sin(\pi s)\int_{-\infty}^{-C}|\lambda|^{-s}\tilde\xi'(\lambda)d\lambda-\hat{\zeta}_2(s)
\end{equation}
can be analitically continued to $\Re s>-1$ with $R(0, \tilde\lambda)=R'_s(0, \tilde\lambda)=0$.
\end{lemma}

For completeness we give a sketch of proof. Using (\ref{trace1}), one has for sufficiently large $\Re s$

$$\zeta(s, \Delta_{\alpha, P}-\tilde\lambda)-\zeta(s, \Delta-\tilde \lambda)=\frac{1}{2\pi i}\int_{c_{\tilde\lambda, \epsilon}}(\lambda-\tilde\lambda)^{-s}{\rm Tr}((\Delta_{\alpha, P}-\lambda)^{-1}-(\Delta-\lambda)^{-1})d\lambda=$$
$$=\int_{c_{\tilde\lambda, \epsilon}}(\lambda-\tilde\lambda)^{-s}\tilde\xi'(\lambda)\,d\lambda=\zeta_1(s)+\zeta_2(s)\,,$$
where
$$\zeta_1(s)=\left\{\int_{-\infty+i\epsilon}^{-C+i\epsilon}-\int_{-\infty-i\epsilon}^{-C-i\epsilon}\right\}(\lambda-\tilde\lambda)^{-s}\tilde\xi'(\lambda)d\lambda\,.$$
It is easy to show (see Lemma 5. 8 in \cite{HK}) that in the limit $\epsilon\to 0$   $\zeta_1(s)$ gives
\begin{equation}\label{eq1}
2i\sin(\pi s)\int_{-\infty}^{-C}|\lambda|^{-s}\tilde\xi'(\lambda)\,d\lambda+2i\sin(\pi s)\int_{-\infty}^{-C}|\lambda|^{-s}\tilde\xi'(\lambda)\rho(s, \tilde \lambda/\lambda)d\lambda\,,
\end{equation}
where $\rho(s, z)=(1+z)^{-s}-1$ and
$$\rho(s, \tilde\lambda/\lambda)=O(|\lambda|^{-1})$$
as $\lambda\to -\infty$.
Using the assumption on the asymptotics of $\tilde\xi(\lambda)$ as $\lambda\to -\infty$ and the obvious relation $\rho(0, z)=0$ one can see that the last term in (\ref{eq1}) can be analytically continued to $\Re s>-1$ and vanishes together with its first derivative w. r. t. $s$ at $s=0$. Denoting it by $R(s, \tilde\lambda)$ one gets the Lemma. $\square$

As it is stated in the introduction the main object we are to study in the present paper is the zeta-regularized determinant of the operator $\Delta_{\alpha, P}-\lambda$.
Let us remind the reader that the usual definition of the zeta-regularized determinant of an operator $A$
\begin{equation}\label{def}{\rm det}\,A=\exp{(-\zeta'(0, A))}\end{equation}
requires analyticity of $\zeta(s, A)$ at $s=0$.

Since the operator zeta-function  $\zeta(s, \Delta-\tilde \lambda)$ is regular at $s=0$ (in fact, it is true in case of $\Delta$ being an arbitrary elliptic differential operator on any compact manifold) and the function $\hat{\zeta}_2(s)$ is  entire,    Lemma \ref{lem} shows that the behavior of the
function $\zeta(s, \Delta_{\alpha, P}-\tilde\lambda)$ at $s=0$ is determined by the properties of the analytic continuation of the term
\begin{equation}\label{TERM} 2i\sin(\pi s)\int_{-\infty}^{-C}|\lambda|^{-s}\tilde\xi'(\lambda)d\lambda\end{equation}
in (\ref{MAIN}). These properties in their turn are determined by the asymptotical behavior of the function $\tilde \xi'(\lambda)$ as $\lambda\to -\infty$.

 It turns out that the latter behavior depends on dimension $d$. In particular, in the next section we will find out that in case $d=2$ the function $\zeta(s, \Delta_{\alpha, P}-\tilde\lambda)$ is not regular at $s=0$, therefore, in order to define ${\rm det}(\Delta_{\alpha, P}-\tilde\lambda)$ one has to use a modified version of (\ref{def}) .

\section{Determinant of pseudo-laplacian on two-dimensional compact manifold}
Let $X$ be a two-dimensional Riemannian manifold, then introducing isothermal local coordinates $(x, y)$  and setting $z=x+iy$, one can write the area element on $X$ as
$$\rho^{-2}(z)|dz|^2$$

The following estimate of the resolvent kernel, $R(z', z; \lambda)$, of the Laplacian on $X$ was found by J. Fay (see \cite{Fay}; Theorem 2.7 on page 38 and the formula preceding Corollary 2.8 on page 39; notice that Fay works with negative Laplacian, so one has to take care of signs when using his formulas).
\begin{lemma}\label{Fay}({\rm J. Fay})
The following equality holds true
\begin{equation}\label{asfay}
-R(z, z'; \lambda)=G_2(r)+\frac{1}{2\pi}\left[\gamma+\log\frac{\sqrt{|\lambda|+1}}{2}\right.\end{equation}
$$\left.-\frac{1}{2(|\lambda|+1)}(1+\frac{4}{3}\rho^{2}(z)\partial^2_{z\bar z}\rho(z))+\hat{R}(z', z; \lambda)\right]\,,
$$
 where $\hat{R}(z', z; \lambda)$ is continuous for $z'$ near $z$, $$\hat {R}(z, z; \lambda)=O(|\lambda|^{-2})$$  uniformly w. r. t. $z\in X$ as $\lambda\to -\infty$;  $r={\rm dist}(z, z')$, $\gamma$ is the Euler constant.
\end{lemma}

Using (\ref{asfay}), we immediately get the following asymptotics of the scattering coefficient $F(\lambda,P)$ as $\lambda\to -\infty$:
\begin{equation}\label{asF2}
F(\lambda, P)=\end{equation}
$$\frac{1}{4\pi}\log(|\lambda|+1)+\frac{\gamma-\log 2}{2\pi}-\frac{1}{4\pi(|\lambda|+1)}\left[1+\frac{4}{3}\rho^{2}(z)\partial^2_{z\bar z}\rho(z)\Big|_{z=z(P)}\right]+O(|\lambda|^{-2})\,.
$$
\begin{remark}{\rm It is easy to check that the expression $\rho^{2}(z)\partial^2_{z\bar z}\rho(z)\Big|_{z=z(P)}$ is independent of the choice of  conformal local parameter $z$ near $P$.}
\end{remark}

Now from (\ref{xsi}) and (\ref{asF2}) it follows that
$$2\pi i\tilde\xi'(\lambda)=-\frac{\frac{1}{4\pi(|\lambda|+1)}-\frac{b}{(|\lambda|+1)^{2}}+O(|\lambda|^{-3})}{\cot \alpha -a-\frac{1}{4\pi}\log(|\lambda|+1)+\frac{b}{|\lambda|+1}+O(|\lambda|^{-2})},$$
where
$a=\frac{1}{2\pi}(\gamma-\log 2)$ and $b=\frac{1}{4\pi}(1+\frac{4}{3}\rho^2\partial^2_{z\bar z}\rho)$.
This implies that for $-\infty<\lambda\leq -C$ one has
\begin{equation}\label{eq3}
2\pi i\tilde\xi'(\lambda)=\frac{1}{|\lambda|(\log|\lambda|-4\pi\cot \alpha+4\pi a)}+f(\lambda)\,,
\end{equation}
with $f(\lambda)=O(|\lambda|^{-2})$ as $\lambda\to -\infty$.
Now knowing (\ref{eq3}), one can study the behaviour of the term (\ref{TERM}) in (\ref{MAIN}).
We have
\begin{equation}\label{TERM1} 2i\sin(\pi s)\int_{-\infty}^{-C}|\lambda|^{-s}\tilde\xi'(\lambda)d\lambda=\end{equation}
$$\frac{\sin(\pi s)}{\pi}\int_{-\infty}^{-C}|\lambda|^{-s-1}\frac{d\lambda}{(\log|\lambda|-4\pi\cot \alpha+4\pi a)}+\frac{\sin(\pi s)}{\pi}\int_{-\infty}^{-C}|\lambda|^{-s}f(\lambda)\,d\lambda\,.
$$

The first integral in the right hand side of (\ref{TERM1}) appeared in (\cite{KLP}, p. 15), where it was observed that it can be easily rewritten through the function
$${\rm Ei}(z)=-\int_{-z}^\infty e^{-y}\frac{dy}{y}=\gamma+\log(-z)+\sum_{k=1}^{\infty}\frac{z^k}{k\cdot k!}\,$$
which leads to the representation
\begin{equation}
\frac{\sin(\pi s)}{\pi}\int_{-\infty}^{-C}|\lambda|^{-s-1}\frac{d\lambda}{(\log|\lambda|-4\pi\cot \alpha+4\pi a)}=
\end{equation}
$$-\frac{\sin(\pi s)}{\pi}e^{-s\kappa}\left[\gamma+\log (s(\log C-\kappa))+ e(s)\right]\,
$$
where $e(s)$ is an entire function such that $e(0)=0$; $\kappa=4\pi\cot\alpha-4\pi a$.
From this we conclude that
\begin{equation}\label{eq5}
\frac{\sin(\pi s)}{\pi}\int_{-\infty}^{-C}|\lambda|^{-s-1}\frac{d\lambda}{(\log|\lambda|-4\pi\cot \alpha+4\pi a)}=-s\log s+g(s)\,\end{equation}
where $g(s)$ is differentiable at $s=0$.

Now (\ref{MAIN}) and (\ref{eq5}) justify the following definition.
\begin{definition} Let $\Delta_{\alpha, P}$ be the pseudo-laplacian on  a two-dimensional compact Riemannian manifold. Then the zeta-regularized determinant of the operator $\Delta_{\alpha, P}-\tilde \lambda$ with $\tilde\lambda\in {\mathbb C}\setminus {\rm Spectrum}(\Delta_{\alpha, P})$
is defined as
\begin{equation}\label{def1}
{\rm det}(\Delta_{\alpha, P}-\tilde \lambda)=\exp\left\{-\frac{d}{ds}\left[\zeta(s, \Delta_{\alpha, P}-\tilde \lambda)+s\log s\right]\Big|_{s=0} \right\}\end{equation}
\end{definition}

We are ready to get our main result: the formula relating ${\rm det}(\Delta_{\alpha, P}-\tilde \lambda)$ to ${\rm det}(\Delta-\tilde\lambda)$.

From (\ref{MAIN}, \ref{TERM}) it follows that
$$\frac{d}{ds}\left[\zeta(s, \Delta_{\alpha, P}-\tilde\lambda) +s\log s-\zeta(s,\Delta-\tilde\lambda)\right]\Big|_{s=0}=$$
$$\frac{d}{ds}\hat{\zeta}_2(s)\Big|_{s=0}+\int_{-\infty}^{-C}f(\lambda)\,d\lambda+$$
$$-\frac{d}{ds}\left\{\frac{\sin \pi s}{\pi}e^{-s\kappa}\left[ \gamma+\log (s(\log C-\kappa))+ e(s)\right]+s\log s                                \right\}\Big|_{s=0}=$$
$$2\pi i\left(\tilde\xi(\tilde\lambda)-\tilde\xi(-C)\right)+\int_{-\infty}^{-C}f(\lambda)\,d\lambda
-\gamma-\log(\log C- \kappa)=
$$
\begin{equation}\label{eq11}
2\pi i \tilde \xi(\tilde \lambda)-\gamma+\end{equation}$$\int_{-\infty}^{-C}f(\lambda)\,d\lambda -   2\pi i\tilde\xi(-C)-\log(\log C-4\pi\cot\alpha+2\gamma-\log 4)\,.
$$
Notice that the expression in the second line of (\ref{eq11}) should not depend on $C$, so one can send  $C$ to $+\infty$ there.
Together with (\ref{asF2}) this gives
\begin{equation}\label{result1}
\frac{d}{ds}\left[\zeta(s, \Delta_{\alpha, P}-\tilde\lambda) +s\log s-\zeta(s,\Delta-\tilde\lambda)\right]\Big|_{s=0}=\end{equation}$$2\pi i \tilde\xi(\tilde\lambda)-\gamma+\log(\sin \alpha/(4\pi))-i\pi\,$$
which implies the comparison formula for the determinants stated in the following theorem.
\begin{theorem}\label{th1}
Let $\tilde\lambda$ do not belong to the union of spectra of $\Delta$ and $\Delta_{\alpha, P}$ and let the zeta-regularized determinant of $\Delta_{\alpha, P}$ be defined as in (\ref{def1}). Then one has the relation
\begin{equation}\label{Mresult2}
{\rm det}(\Delta_{\alpha, P}-\tilde\lambda)= -4\pi e^{\gamma}(\cot \alpha- F(\tilde\lambda, P)){\rm det}(\Delta-\tilde\lambda)\, .\end{equation}
\end{theorem}

Observe now that $0$ is the simple eigenvalue of $\Delta$ and, therefore, it follows from Theorem 2 in \cite{YCdV}  that $0$ does not belong to the spectrum of the operator  $\Delta_{\alpha, P}$ and that   $\Delta_{\alpha, P}$ has one strictly negative simple eigenvalue.
Thus, the determinant in the left hand side of (\ref{Mresult2}) is well defined for $\tilde\lambda=0$, whereas
the determinant at the right hand side has the asymtotics
\begin{equation}\label{000}
{\rm det}(\Delta-\tilde\lambda)\sim (-\tilde \lambda){\rm det}^*\Delta\, \end{equation}
as $\tilde\lambda\to 0-$. Here ${\rm det}^*\Delta$ is the modified determinant of an operator with zero mode.

From the standard asymptotics
$$-R(x, y; \lambda)=\frac{1}{{\rm Vol}(X)}\frac{1}{\lambda}+G_2(r)+O(1)$$
as $\lambda\to 0$ and $x\to y$ one gets the asymptotics
\begin{equation}\label{1111}
F(\lambda, P)=\frac{1}{{\rm Vol}(X)}\frac{1}{\lambda}+O(1)
\end{equation}
as $\lambda\to 0$.
Now sending $\tilde\lambda\to 0-$ in (\ref{Mresult2}) and using \ref{000} and \ref{1111}
we get the following corollary of the Theorem \ref{th1}.
\begin{corollary} The following relation holds true
\begin{equation}
{\rm det}\Delta_{\alpha, P}= -\frac{4\pi e^{\gamma}}{{\rm Vol}(X)}{\rm det}^*\Delta\,.
\end{equation}

\end{corollary}

\section{Determinant of pseudo-laplacian on three-dimensional manifolds}
 Let $X$ be a three-dimensional compact Riemannian manifold.  We start with the Lemma describing the asymptotical behavior of the scattering coefficient as $\lambda \to -\infty$.

\begin{lemma}\label{GlLemm} One has the asymptotics
\begin{equation}\label{asy} F(\lambda; P)=\frac{1}{4\pi}\sqrt{-\lambda}+c_1(P)\frac{1}{\sqrt{-\lambda}}+O(|\lambda|^{-1})\end{equation}
as $\lambda\to -\infty$
\end{lemma}

{\bf Proof.}
Consider Minakshisundaram-Pleijel asymptotic expansion (\cite{MP})
\begin{equation}\label{MP1} H(x, P; t)=(4\pi t)^{-3/2}e^{-d(x, P)^2/(4t)}\sum_{k=0}^\infty u_k(x, P)t^k\end{equation}
for the heat kernel in a small vicinity of $P$, here $d(x, P)$ is the geodesic distance from $x$ to $P$, functions $u_k( \cdot , P)$ are smooth in a vicinity of P, the equality is understood in the sense of asymptotic expansions.
We will make use of the standard relation
\begin{equation}\label{stR}R(x, y; \lambda)=\int_0^{+\infty}H(x, y; t)e^{\lambda t}\,dt\,.\end{equation}
Let us first truncate the sum (\ref{MP1}) at some fixed $k=N+1$  so that the remainder, $r_n$,  is $O(t^{N})$.
Defining
\[
\tilde{R}_N(x,P; -\lambda) := \int_0^\infty r_n(t,x,P) e^{t\lambda} dt\,,
\]
we see that
$$\tilde{R}_N(x,P; \lambda)=O(|\lambda|^{-(N+1)})$$
as $\lambda\to -\infty$ uniformly w. r. t. $x$ belonging to a small vicinity of $P$.

Now, for each $0\leq k\leq N+1$ we have to address the following quantity

\[
R_k(x,P; \lambda) := \frac{u_k(x,y)}{(4\pi)^{3/2}} \int_0^\infty t^{k-\frac{3}{2}}e^{-\frac{d(x, P)^2}{4t}}e^{\lambda t} dt.
\]

According to identity (\ref{Ba}) below one has
$$R_0(x, P; \lambda)=\frac{u_0(x, P)}{(4\pi)^{3/2}}\frac{2\sqrt{\pi}}{d(x, P)}e^{-d(x, P)\sqrt{-\lambda}}=$$
\begin{equation}\label{te1}
\frac{1}{4\pi d(x, P)}-\frac{1}{4\pi}\sqrt{-\lambda}+o(1),
\end{equation}
as $d(x, P)\to 0$.
For $k\geq 1$ one has
$$R_k(x, P; \lambda)=\frac{u_k(x, P)}{(4\pi)^{3/2}}2^{3/2-k}\left(\frac{d(x, P)}{\sqrt{-\lambda}} \right)^{k-1/2}K_{k-\frac{1}{2}}(d(x, P)\sqrt{-\lambda})=$$
\begin{equation}\label{te2}
-c_k(P)\frac{1}{(\sqrt{-\lambda})^{2k-1}}+o(1)\end{equation}
as $d(x, P)\to 0$ (see \cite{Ba}, p. 146, f-la 29). Now (\ref{asy}) follows from (\ref{stR}), (\ref{te1}) and (\ref{te2}).
 $\square$

Now from Lemma \ref{GlLemm} it follows that
\begin{equation}\label{glavn}
2\pi i\tilde\xi'(\lambda)=-\frac{1}{2\lambda}+O(|\lambda|^{-3/2})\end{equation}
as $\lambda\to -\infty$, therefore, one can rewrite (\ref{TERM}) as

\begin{equation}\label{T2}\frac{\sin(\pi s)}{\pi}\left\{\int_{-\infty}^{-C}|\lambda|^{-s}(2\pi i \tilde\xi'(\lambda)+\frac{1}{2\lambda})d\lambda+\frac{C^{-s}}{2s}\right\}\end{equation}
which is obviously analytic in $\Re s>-\frac{1}{2}$.
Thus, it follows from (\ref{MAIN}) that the function $\zeta(s, \Delta_{\alpha, P}-\tilde\lambda)$ is regular at $s=0$ and one can introduce the usual zeta-regularization
$${\rm det}(\Delta_{\alpha, P}-\tilde\lambda)=\exp\{-\zeta'(0, \Delta_{\alpha, P}-\tilde\lambda)\}$$
of ${\rm det}(\Delta_{\alpha, P}-\tilde\lambda)$.

Moreover, differentiating (\ref{MAIN}) with respect to $s$ at $s=0$ similarly to (\ref{eq11}) we get
$$\frac{d}{ds}\left[\zeta(s, \Delta_{\alpha, P}-\tilde\lambda) -\zeta(s,\Delta-\tilde\lambda)\right]\Big|_{s=0}=$$
$$2\pi i(\tilde\xi(\tilde\lambda)-\tilde\xi(-C))+\int_{-\infty}^{-C}(2\pi i \tilde\xi'(\lambda)+\frac{1}{2\lambda})d\lambda-\frac{1}{2}\log C=$$
which reduces after sending $-C\to-\infty$ to
$$2\pi i\tilde\xi(\tilde\lambda)+\log\sin \alpha-\log (4\pi)+i\pi=-\log(\cot\alpha-F(\lambda; P))-\log(4\pi)+i\pi\,$$
which implies the following theorem.
\begin{theorem} Let $\Delta_{\alpha, P}$ be the pseudo-laplacian on $X$ and $\tilde\lambda\in {\mathbb C}\setminus ({\rm Spectrum}(\Delta)\cup
{\rm Spectrum}(\Delta_{\alpha, P}))$. Then
\begin{equation}\label{SphereDet}
{\rm det}(\Delta_{\alpha, P}-\tilde\lambda)=-4\pi(\cot\alpha-F(\tilde\lambda; P)){\rm det}(\Delta-\tilde\lambda)\,.
\end{equation}
\end{theorem}

Sending $\tilde\lambda\to 0$ and noticing that relation ({\ref {1111}) holds also in case $d=3$ we get the following corollary.
\begin{corollary}
\begin{equation}
{\rm det}\Delta_{\alpha, P}=-\frac{4\pi}{{\rm Vol}(X)}{\rm det}^*\Delta\,.
\end{equation}

\end{corollary}

  In what follows we consider two examples of three-dimensional compact Riemannian manifolds for which there exist explicit expressions for the resolvent kernels: a flat torus and the round (unit) $3d$-sphere. These manifolds are homogeneous, so, as it is shown in \cite{YCdV}, the scattering coefficient $F(\lambda, P)$ is $P$-independent.

{\bf Example 1: Round $3d$-sphere.}
\begin{lemma}\label{sph} Let $X=S^3$ with usual round metric.
  Then there is the following explicit expression for scattering coefficient
  \begin{equation}\label{Tayeb}F(\lambda)=\frac{1}{4\pi}\coth\left(\pi\sqrt{-\lambda-1}\,\right)\cdot\sqrt{-\lambda-1}\end{equation}
and, therefore, one has the following asymptotics as $\lambda\to -\infty$
\begin{equation}\label{sph1}
F(\lambda)=\frac{1}{4\pi}\sqrt{|\lambda|-1}+O(|\lambda|^{-\infty})\,.
\end{equation}
\begin{remark}{\rm The possibility of finding an explicit expression for $F(\lambda)$ for $S^3$ was mentioned in \cite{YCdV}. However we failed to find (\ref{Tayeb}) in the literature.}
\end{remark}
\end{lemma}
{\bf Proof.} We will make use the well-known identity (see, e. g., \cite{Ba}, p. 146, f-la 28):
\begin{equation}\label{Ba}
\int_{0}^{+\infty}e^{\lambda t}t^{-3/2}e^{-\frac{d^2}{4t}}\,dt=2\frac{\sqrt{\pi}}{|d|}e^{-|d|\sqrt{-\lambda}};
\end{equation}
for $\lambda<0$ and $d\in{\mathbb R}$ and the following explicit formula for the operator kernel $e^{-t}H(x, y; t)$ of the operator
$e^{-t(\Delta+1)}$, where $\Delta$ is the (positive) Laplacian on $S^3$ (see \cite{CT}, (2.29)):
\begin{equation}\label{Taylor}
e^{-t}H(x, y; t)=-\frac{1}{2\pi}\frac{1}{\sin d(x, y)}\frac{\partial}{\partial z}\Big|_{z=d(x,y)}\Theta(z, t)\,.
\end{equation}
Here $d(x, y)$ is the geodesic distance between $x, y\in S^3$ and
$$\Theta(z, t)=\frac{1}{\sqrt{4\pi t}}\sum_{k=-\infty}^{+\infty}e^{-(z+2k\pi)^2/4t}$$
is the theta-function.

Denoting $d(x, y)$ by $\theta$ and using (\ref{Taylor}) and (\ref{Ba}), one gets
$$R(x, y; \lambda-1)=\int_{0}^{+\infty}e^{\lambda t}e^{-t}H(x, y; t)\,dt=$$$$\frac{1}{4\pi}\frac{1}{\sin \theta}\left(-\sum_{k<0}e^{(\theta+2k\pi)\sqrt{-\lambda}}+\sum_{k\geq 0}e^{-(\theta+2k\pi)\sqrt{-\lambda}}\right)=$$
$$\frac{1}{4\pi}\frac{1}{\sin \theta}\frac{1}{1-e^{-2\pi\sqrt{-\lambda}}}\left[-e^{-2\pi\sqrt{-\lambda}}e^{\theta\sqrt{-\lambda}}+e^{-\theta\sqrt{-\lambda}}
\right]=$$
\begin{equation}
\frac{1}{4\pi\theta}-\frac{1}{4\pi}\frac{1+e^{-2\pi\sqrt{-\lambda}}}{1-e^{-2\pi\sqrt{-\lambda}}}\sqrt{-\lambda}+o(1)
\end{equation}
as $\theta\to 0$, which implies the Lemma.  $\square$

{\bf Example 2: Flat $3d$-tori.}Let $\{{\bf A, B, C}\}$ be a basis of ${\mathbb R}^3$ and let $T^3$ be the quotient of ${\mathbb R}^3$ by the lattice $\{m{\bf A}+n{\bf B}+l{\bf C}: (m, n, l)\in {\mathbb Z}^3\}$ provided with the usual flat metric.

Notice that the free resolvent kernel
in $R^3$ is $$\frac{e^{-\sqrt{-\lambda}||x-y||}}{4\pi||x-y||}$$ and, therefore,
\begin{equation}\label{exactres}
R(x, y; \lambda)=\frac{e^{-\sqrt{-\lambda}||x-y||}}{4\pi||x-y||}+\frac{1}{4\pi}
\sum_{(m, n, l)\in {\mathbb Z}^3\setminus(0,0,0)}\frac{e^{-\sqrt{-\lambda}||x-y+m{\bf A}+n{\bf B}+l{\bf C}||}}{||x-y+m{\bf A}+n{\bf B}+l{\bf C}||}\,.
\end{equation}
From (\ref{exactres})  it follows that
$$F(\lambda)= \frac{1}{4\pi}\sqrt{-\lambda}-\frac{1}{4\pi}\sum_{(m, n, l)\in {\mathbb Z}^3\setminus(0,0,0)}\frac{e^{-\sqrt{-\lambda}||m{\bf A}+n{\bf B}+l{\bf C}||}}{||m{\bf A}+n{\bf B}+l{\bf C}||}=$$$$\frac{1}{4\pi}\sqrt{-\lambda}+O(|\lambda|^{-\infty})$$
as $\lambda\to-\infty$.

\begin{remark}{\rm It should be noted that explicit expressions for ${\rm det}^*\Delta$ in case $X=S^3$ and $X=T^3$ are given in
\cite{Kum} and \cite{FG2003}.}
\end{remark}

\end{document}